\documentclass[a4paper]{amsart}

\usepackage{amsmath}
\usepackage{amssymb}
\usepackage{amsfonts}
\usepackage{amsthm}
\usepackage[top=28truemm,bottom=30truemm,left=29truemm,right=29truemm]{geometry}
\usepackage[fontsize=11pt]{scrextend}
\usepackage{comment}
\usepackage{bm}
\usepackage{mathrsfs}

\newtheorem{thm}{Theorem}[section]
\newtheorem{prop}{Proposition}[section]

\newtheorem{cor}{Corollary}[section]

\def\supp{{\rm supp}\,}

\def\ls{\lesssim}
\def\bs{\backslash}

\def\dr{\partial_r}
\def\dt{\partial_t}
\def\ds{\partial_s}

\def\Com{\mathbb C}
\def\R{\mathbb R}
\def\N{\mathbb N}
\def\Z{\mathbb Z}

\def\D{\mathcal{D}}
\def\E{\mathcal{E}}
\def\G{\mathcal{G}}
\def\H{\mathcal{H}}
\def\M{\mathcal{M}}

\def\P{\mathcal{P}}
\def\Q{\mathcal{Q}}

\def\cN{\mathcal{N}}
\def\L{\mathcal{L}}
\def\S{\mathcal{S}}

\def\U{\mathcal{U}}
\def\V{\mathcal{V}}

\makeatletter
\@addtoreset{equation}{section}
\makeatother

\title[The 2D stationary Navier-Stokes equations]{Existence of the stationary Navier-Stokes flow in $\R^2$
around a radial flow}

\author[Yasunori Maekawa \and Hiroyuki Tsurumi]{Yasunori Maekawa$^1$ \and Hiroyuki Tsurumi$^2$}

\thanks{$^{1,2}$Department of Mathematics, Graduate School of Science, Kyoto University, Kyoto 606--8502, Japan.
\\\indent
$^1$maekawa.yasunori.3n@kyoto-u.ac.jp 
\\\indent
$^2$tsurumi.hiroyuki.23h@st.kyoto-u.ac.jp}

\begin{document}
\maketitle

\begin{abstract}
We consider the stationary Navier-Stokes equations on the whole plane $\R^2$.
We show that for a given small and smooth external force around a radial flow, there exists a classical solution decaying like $|x|^{-1}$.
In our result, it is not necessary to impose any symmetric conditions on external forces.
\end{abstract}

\section{Introduction}

We consider the two-dimensional stationary Navier-Stokes equations
\begin{equation}
\label{NS}
\left\{
  \begin{array}{ll}
    -\Delta u+u\cdot\nabla u+\nabla p = F,\hspace{5pt}  &x\in \R^2, \\
    \nabla\cdot u=0 \hspace{5pt}  &x\in \R^2.
  \end{array}
\right.
\end{equation}
Here $u=(u_1(x),u_2(x)) $ and $p=p(x)$ denote the unknown velocity vector and the unknown pressure of the fluid at the point $x=(x_1,x_2)\in \R^2$, respectively, while $F=(F_1(x),F_2(x))$ is the given external force. Here $u\cdot\nabla u := \sum_{j=1,2}u_j\partial_{x_j}u$.

For three or higher dimension case, the system \eqref{NS} have been studied well.
For instance, Leray \cite{Le1933} and Ladyzhenskaya \cite{La1959} showed the existence of strong solutions to \eqref{NS}, and Heywood \cite{He1970} constructed solutions of \eqref{NS} as a limit of solutions of the non-stationary Navier-Stokes equations.
Later on, various researchers have found scaling invariant spaces of $F$ guaranteeing the existence of solutions in $\R^3$, such as Chen \cite{Ch1993} in the Lebesgue space, Kozono-Yamazaki \cite{KY1995} in the Morrey space, and Kaneko-Kozono-Shimizu \cite{KKS2019} in the Besov space.
These theories, however, cannot be applied for the two dimension case,
since it is hard to estimate the advection term $u\cdot \nabla u$ in scaling invariant function spaces.
Hence, until now, the two dimension case has been considered independently of higher dimension ones.

Actually, in two dimensional exterior domains, the system \eqref{NS} has been investigated for many years. 
Classically, Finn-Smith \cite{FS1967} considered a solution $u$ to \eqref{NS} with $F=0$, which satisfies $u(x)\to v$ as $|x|\to\infty$ with some constant vector $v\neq 0$.
Then Amick \cite{Am1984} showed the existence of solutions for given external forces when an exterior domain is invariant under the transformation $(x_1, x_2)\mapsto (-x_1, x_2)$, and this work was later generalized  by Pileckas-Russo \cite{PR2012} for example.
For such solutions, Galdi-Sohr \cite{GS1995} and \cite{PR2012} precisely investigated the regularity and asymptotic behavior under some symmetric conditions on $F$.
After that, Yamazaki \cite{Ya2016} showed the existence of strong and weak solutions for every $F=\nabla\cdot \tilde F$ with $\tilde F\in (L^2(\Omega))^4$ on an exterior domain $\Omega$ which is invariant under some symmetric group action.
Especially in the exterior of the disk, Hillairet-Wittwer \cite{HW2013} constructed stationary solutions around the radial and rotational flow $\mu x^\perp/|x|^2$ with a sufficiently large constant $\mu$. 
The  flow $\mu x^\perp/|x|^2$ is the well-known exact solution decaying in the scale-critical order $O(|x|^{-1})$ under the zero flux condition, and the key obserbation in \cite{HW2013} is that the vorticity transport by the flow $\mu x^\perp/|x|^2$ leads to the essential change of the decay structure for the fundamental solution of its linearization. 
The largeness condition of $\mu$ plays a key role here, and in particular, it is not known in general whether or not the stationary solution exists when $\mu$ is not large. 
We note that the stability of $\mu x^\perp/|x|^2$ for small perturbations was proved by M. \cite{Ma2017} when $|\mu|$ is small enough; see also Higaki \cite{Hi2019} for the genelarization of this stability result. 
Recently, Higaki-M.-Nakahara \cite{HMN2017, HMN2018} showed the existence and uniqueness of stationary solutions around a rotating obstacle with small external forces and low speed of rotation,
and later on, Gallagher-Higaki-M. \cite{GHM2019} generalized it with high speed.

On the other hand, on the whole plane $\mathbb{R}^2$, Yamazaki \cite{Ya2009} showed the existence of small solutions to \eqref{NS} in the weak $L^2$ space when given external forces decay sufficiently and have some symmetric properties, such as $f_1(x_2,x_1)=f_2(x_1,x_2)$ for example.
Galdi-Yamazaki \cite{GY2015} later showed the stability of the above solutions more precisely.
For its uniqueness, Nakatsuka \cite{Na2015} constructed more general theory.
Recently, Guillod \cite{Gu2015} showed the existence of a pair $(u, F)$ solving \eqref{NS}, where $F$ is dependent on $u$ and is constructed around an arbitrarily given small function $k$ having zero integral and decaying faster than $|x|^{-3}$.
Moreover, Guillod-Wittwer \cite{GW2015} found solutions to \eqref{NS} which is scaling invariant with respect to a rotation conversion.

This paper contributes to the open problem related to \cite{HW2013}, in the sense that at least in the $\R^2$ setting, rather than the exterior problem, we show that the stationary solution behaves like $\mu x^\perp/|x|^2$ for $|x|\gg 1$ is constructed for small but nonzero $\mu$ for a suitable class of the external forces. 
Let us suppose that a external force $F$ is divergence-free and can be expressed as $F=\nabla^\perp \phi:=(\partial_{x_2}\phi, -\partial_{x_1}\phi)$ with some flow $\phi$.
Actually, by the Helmholtz decomposition, we can write $F$ as a sum of the rotation-free term and divergence-free one.
Then the curl-free term is absorbed in the pressure term $\nabla p$, and the divergence-free one is written as $\nabla^\perp \phi$ if $F$ decays sufficiently.
Under such a situation, we will show that for every small compact supported radial flow $\phi_*$ having non-zero integral and its smaller perturbation $\varphi$ decaying faster than $|x|^{-2}$, there exist solutions $(u, p)$ of \eqref{NS} for $F=\nabla^\perp(\phi_*+\varphi)$. 
This decay condition implies
\begin{equation*}
\int_{\R^2}F(x)dx=0,
\end{equation*}
which means that the source term $F$ provides no net force.
Here we should note that Guillod-Wittwer \cite{GW2015A} showed more complicated behaivior of a solution in the far-field when $F$ provides non-zero net force.
In our result, we also impose the smoothness of $\phi_*$ and $\varphi$, so that we may obtain $u$ as a classical solution in $C^2(\R^2; \R)$.

The standard approach to \eqref{NS} is to analyze the vorticity-streamfunction system as below, which is equivalent with \eqref{NS} in a suitable functional framework.
\begin{equation}
\label{VS}
\left\{
  \begin{array}{ll}
    \Delta\psi= -\omega, \hspace{5pt}  &x\in \R^2, \\
    \Delta\omega=\nabla\times(\nabla^\perp \psi\cdot \nabla(\nabla^\perp \psi))
+\Delta \phi, \hspace{5pt}  &x\in \R^2.
  \end{array}
\right.
\end{equation}
Here $\psi=\psi(x)$ is the stream function, which generates the divergence-free flow as $u=\nabla^\perp\psi$, while $\omega:=\nabla\times u:=\partial_{x_1}u_2-\partial_{x_2}u_1$ is the vorticity field. 
In addition, we convert orthogonal coordinates $(x_1, x_2)$ to polar ones $(r,\theta)$, and consider the Fourier series $\psi(r,\theta)=\sum_{n\in\Z} \psi_n(r) e^{in\theta}$ and $\omega(r,\theta)= \sum_{n\in\Z} \omega_n(r) e^{in\theta}$ with respect to the angular valuable $\theta$.
Then \eqref{VS} is expressed as the system of ordinary differential equations of $\{\psi_n\}_{n\in\Z}$ and $\{\omega_n\}_{n\in\Z}$ with the radius valuable $r$. 
For this method, we especially refer to Hillairet-Wittwer \cite{HW2013}.

This paper is organized as follows. 
In the next section, we will define the Fourier modes of vector fields on polar coordinates and some important function spaces of them, and we will state our main theorem. 
In the third section, we will first state the policy of the proof, and then we show some required propositions.

\section{Main result}

We will analyze the equations in the polar coordinates. 
Actually, some key structures are found by decomposing the system \eqref{VS} into the Fourier mode with respect to the angular valuable.
For later use, let us introduce the Fourier series, as
\begin{equation*}
f(r\cos\theta, r\sin\theta) = \sum_{n\in\Z}f_n(r)e^{in\theta},
\hspace{10pt} (r,\theta)\in [0,\infty)\times[0,2\pi),
\end{equation*}
where $f_n$ denotes the $n$-mode of $f$ defined by
\begin{equation*}
f_n(r)=\frac1{2\pi} \int_0^{2\pi} f(r\cos\theta, r\sin\theta)e^{-in\theta}d\theta,
\hspace{10pt} 0\leq r<\infty.
\end{equation*}
In what follows, we write $\hat f=(f_n)_{n\in\Z}$ and $\tilde f=(f_n)_{n\in \Z\backslash\{0\}}$ for the Fourier mode of $f$. 
We now introduce the following function spaces for the Fourier mode.
For simplicity, we write continuous and smooth function spaces as
\begin{equation*}
C^m:= C^m([0,\infty);\Com),\hspace{5pt} \hat C^m:=(C^m)^\Z,\hspace{5pt}
 \tilde C^m:=(C^m)^{\Z\backslash\{0\}}
\end{equation*}
for $m\in\N\cup\{0\}$.
Let $\alpha>0$ and $\kappa>1$. 
For a function $f\in C^l$, define the weight function $\M ^l_{n; \alpha,\kappa}[f]$ as
\begin{equation*}
\M ^l_{n; \alpha, \kappa}[f](r)
:=(1+r)^{\alpha+l}(1+|n|)^{\kappa-l}|\partial^l_r f(r)|,\hspace{10pt}n\in\Z,\  l\in\N\cup\{0\},\ l<\kappa.
\end{equation*}
We note here that $\alpha$ counts the decay in $r$ and $\kappa$ the decay in $n$. In addition, the order $l$ of derivative makes the decay in $r$ faster and the decay in $n$ slower.
In this paper, we applied this weight $\M ^l_{n; \alpha,\kappa}$ for the Fourier mode $f_n$ with the same $n\in\Z$, so that in what follows, we use the abbreviation as
\begin{equation*}
\M ^l_{\alpha, \kappa}[f_n](r):=\M ^l_{n; \alpha, \kappa}[f_n](r).
\end{equation*}
In association with this weight function, we define some norms and spaces as follows.
As for spaces of the vorticity and external force, we define the norm
\begin{equation*}
\|\hat f\|_{{\U}^m_{\alpha,\kappa}}
:=\sum_{l=0}^m \sup_{n\in\Z} \sup_{r\geq0} \M ^l_{\alpha, \kappa}[f_n](r),\hspace{10pt} \hat f\in \hat C^m,
\ m\in\N\cup\{0\},\ m<\kappa
\end{equation*}
and set
\begin{equation*}
\U^m_{\alpha,\kappa}:= \left\{\hat f\in\hat C^m; \hspace{5pt}
\|\hat f\|_{\U^m_{\alpha,\kappa}}<\infty,\hspace{5pt} f_{-n}= \overline{f_n}\ \forall n\in\Z,
\hspace{5pt}  f_{n'}(0)=0\ \forall n'\in\Z\bs\{0\}\right\}.
\end{equation*}
Here $ \overline{f}$ denotes the complex conjugate of $f$.
Actually, for each component $\hat f\in \U^m_{\alpha,\kappa}$, we see for large $r$ that
\begin{equation*}
|\dr^l f_n(r)|\leq C(1+r)^{-\alpha-l},\hspace{10pt}\forall n\in \Z,\ 0\leq l\leq m,
\end{equation*}
while for small $r$ that
\begin{equation*}
|f_n(r)|\leq Cr,\hspace{10pt}\forall n\in \Z\bs\{0\}.
\end{equation*}
On the other hand, as a space of stream functions, we set 
\begin{equation*}
\V^m_{\alpha,\kappa}:=
\left\{\hat f=(f_0, \tilde f)\in \V^m_0\times \tilde{\V}^m_{\alpha,\kappa}; \hspace{5pt} 
\|\hat f\|_{\U^m_{\alpha,\kappa}}:=\|f_0\|_{\V^m_0}+\|\tilde f\|_{\tilde{\V}^m_{\alpha,\kappa}}<\infty
\right\},
\end{equation*}
where
\begin{equation*}
\begin{split}
&\V^m_0
:= \left\{
f_0\in C^m; \hspace{5pt}
\|f_0\|_{\V^m_0}:=\sum_{l=1}^m \sup_{r\geq0} \M ^l_{0, 0}[f_0](r)<\infty,\ f_0=\overline{f_0}
\right\},\\
& \tilde{\V}^m_{\alpha,\kappa}
:=\left\{\tilde f\in \tilde C^m; \hspace{5pt}
\|\tilde f\|_{\tilde{\V}^m_{\alpha,\kappa}}:=\sum_{l=0}^m \sup_{n\in\Z\backslash\{0\}}\sup_{r\geq0} \M ^l_{\alpha, \kappa}[f_n](r)<\infty,\ f_{-n}= \overline{f_n}\ \forall n\in \Z\right\}.
\end{split}
\end{equation*}
We note that as for $\hat f\in \V^m_{\alpha,\kappa}$, we do not consider the zero mode $f_0$ itself, but consider its differential. 
In addition, each differential $\dr^l f_0$ decays as $|\dr^l f_0(r)|\ls (1+r)^{-l}$, while that of non-zero mode $\dr^l f_n$ decays as $|\dr^l f_n(r)|\ls (1+r)^{-\alpha-l}$.

We easily see that the above $\U^m_{\alpha,\kappa}$ is a Banach space with its norm.
Moreover, there hold the embeddings $\U^m_{\alpha,\kappa}\subset \U^{m'}_{\alpha',\kappa'}$ and 
$\V^m_{\alpha,\kappa}\subset \V^{m'}_{\alpha',\kappa'}$ for every $0<\alpha'\leq \alpha$, $1<\kappa'\leq \kappa$, and $1\leq m'\leq m$ with the estimates
\begin{equation*}
\|\hat f\|_{\U^{m'}_{\alpha',\kappa'}}\leq \|\hat f\|_{\U^m_{\alpha,\kappa}},\hspace{5pt}
\| \hat f \|_{\V^{m'}_{\alpha',\kappa'}}\leq \|\hat f\|_{\V^m_{\alpha,\kappa}}.
\end{equation*}

Our main result now reads:
\begin{thm}
\label{thm:1}
There exists a constant $\delta>0$ such that the following statement holds.

Let $\phi_*\in C^2_c(\R^2;\R)$ be a radial function expressed as $\phi_*(r)\in C_c^2([0,\infty);\R)$ in polar coordinates $(r,\theta)$, and let $R_*\geq 1$ be such that $\supp(\phi_*)\subset [0,R_*]$. 
Moreover, define
\begin{equation*}
\mu_*:=\int_0^{R_*} s\phi_*(s)ds, \hspace{10pt}
\rho_*:=\sqrt2 \left[ \left\{ 1+ \left(\frac{\mu_*}2\right)^2 \right\}^{\frac12}+1 \right]^\frac12 -2,
\end{equation*}
\begin{equation*}
\nu_*:=\sup_{0\leq r\leq R_*}\left|\int_0^r s\phi_*(s)ds \right|+
\sup_{0\leq r\leq R_*}(1+r)^2|\phi_*(r)|+\sup_{0\leq r\leq R_*}(1+r)^2|\dr\phi_*(r)|,
\end{equation*}
and suppose that
\begin{equation}
\label{condition}
\mu_*\neq0,\hspace{10pt} R_*^{\rho_*}\nu_*<\delta.
\end{equation}
Then for every $0<\alpha<\min \{1/2, \rho_*\}$, there exists a constant $\varepsilon=\varepsilon(R_*^{\rho_*}, \alpha)>0$ such that 
for every $\varphi \in C^2(\R^2;\R)$ whose Fourier mode satisfies
$\|\hat\varphi\|_{\U^1_{\alpha+2,\kappa+2}}<\varepsilon$ for some $\kappa>1$, 
there exists a solution $u\in C^2(\R^2;\R^2)$ of \eqref{NS} for an external force $F=\nabla^\perp(\phi_*+\varphi)$ having the decay property
\begin{equation}
\label{decay}
\sup_{r\geq0}\sup_{0\leq\theta<2\pi} (1+r)^{1+\alpha}\left|u(r,\theta)-\left(\frac1r\int_0^r s(\phi_*+\varphi_0)(s)ds \right) e_\theta\right|<\infty,
\end{equation}
where $\varphi_0$ denotes the Fourier zero-mode of $\varphi$, and $e_\theta:=(-\sin\theta, \cos\theta)$ denotes the basis for the direction of increasing angle in polar coordinates.
\end{thm}

Here we remark that $\mu_*$ is also written as $(2\pi)^{-1}\int_{\R^2}\phi_*(x)dx$.
On the other hand, $\rho_*=\Re[(4+2i\mu_*)^{1/2}]-2$ ($\Re[z]$ denotes the real part of $z$) appears in the partial linearization of the system regarding to the Fourier $\pm 2$ modes.
The asymptotic estimate \eqref{decay} implies that the solution $u$ behaves like the radial and rotational flow $c x^\perp/|x|^2$ with some $c\in \R$ that decays in the scale-critical order $O(|x|^{-1})$. 
From the proof we know that the constant $c$ is small but nonzero. 
Hence our result is contrasting to \cite{HW2013}, where the constant $c$ must be large enough in the case of the exterior disk. 
The key observation in the proof of Theorem \ref{thm:1} is that the system \eqref{VS} has two opposite aspects; the one is related to the analysis in the Fourier $\pm 2$ modes with respect to the angular variable $\theta$, where one needs to use the effect of the vorticity transport by the flow $\mu_* x^\perp/|x|^2$ as in \cite{HW2013} to avoid the appearance of the logarithmic loss from the scale-critical decay pointed out by Guillod \cite[Section 3]{Gu2015}. 
The other is related to the Fourier $\pm 1$ modes, where we find the key cancellation property in the nonlinear term $\nabla \times \, (u\cdot \nabla u)$ that seems to be avairable only by regarding the linearized term around  $\mu_* x^\perp/|x|^2$ as the perturbation. 
Hence, we build up the iteration scheme by taking this observation into account, that is, the transport term by the flow $\mu_* x^\perp/|x|^2$ is incorporated as the principal term for the Fourier $\pm 2$ modes, while this term is handled as the perturbation for the other Fourier modes and we use the smallness of $\mu_*$ and the cancellation property in the Fourier $\pm 1$ modes.
The smallness condition of \eqref{condition} is then needed to close the linear estimate.
Another advantage of our result is that there is no restriction on $\varphi$ regarding to its structure such as symmetry, while the previous studies \cite{GY2015, Na2015, Ya2009, Ya2016} require such structural assumptions.


\section{Proof}

\subsection{Outline}

We now fix $\phi_*\in C_c^2([0,\infty);\R)$ satisfying $\supp (\phi_*)\subset[0,R_*]$ and $\mu_*\neq 0$.
Since $\rho_*<\mu_*<\nu_*$, and since the smallness condition in \eqref{condition} should be satisfied, we assume that $0<\rho_*<1$ in what follows.

In terms of the polar coordinates, the vorticity-streamfunction system \eqref{VS} is expressed as
\begin{equation}
\label{VS2}
\left\{
  \begin{array}{ll}
  \displaystyle  \Delta_{r,\theta}\psi= -\omega, \hspace{5pt}  &(r,\theta)\in[0,\infty)\times [0,2\pi), \\
 \displaystyle   \Delta_{r,\theta}\omega=
G+\Delta_{r,\theta} \phi, \hspace{5pt}  &(r,\theta)\in[0,\infty)\times [0,2\pi).
  \end{array}
\right.
\end{equation}
Here $\Delta_{r,\theta}:=\partial^2_r+(1/r)\dr+(1/r^2)\partial^2_\theta$ denotes the Laplacian in the polar coordinates, and
\begin{equation}
\label{G1}
G
:=-\frac{1}{r^2}\dr(r\dr D+ \partial_\theta E) +\frac{1}{r^3}(1+\partial^2_{\theta})D,
\end{equation}
where
\begin{equation*}
D:=\dr\psi\partial_\theta\psi,\hspace{10pt}
E:=\frac{(\partial_\theta\psi)^2}{r}-r(\dr\psi)^2.
\end{equation*}
On the other hand, since
\begin{equation*}
\nabla\times(\nabla^\perp \psi\cdot \nabla(\nabla^\perp \psi))=\nabla^\perp \psi\cdot\nabla\omega
\end{equation*}
by $\nabla\cdot \nabla^\perp \psi=0$, we also see that
\begin{equation}
\label{G2}
G=\frac1r (\partial_\theta \psi\dr\omega- \dr\psi\partial_\theta \omega).
\end{equation}

Since $\phi_*$ depends only on $r$, we see that
\begin{equation*}
(\psi, \omega)=(\psi_*(r), \omega_*(r)):=\left( -\int_0^r \frac1s \int_0^s t\phi_*(t)dtds,\  \phi_*(r) \right)
\end{equation*}
are exact classical solutions of \eqref{VS2} for $\phi=\phi_*$.
Therefore, for given small $\varphi\in C^2(\R^2;\R)$, we aim to construct solutions of \eqref{VS2} for $\phi=\phi_*+\varphi$ such as
\begin{equation}
\label{solution}
\left\{
\begin{array}{ll}
&\displaystyle \psi(r,\theta)=\psi_*(r)+\gamma(r,\theta)=\psi_*(r)+\sum_{n\in\Z}\gamma_n(r)e^{in\theta}, \\
&\displaystyle \omega(r,\theta)=\omega_*(r)+w(r,\theta)=\omega_*(r)+\sum_{n\in\Z}w_n(r)e^{in\theta},
\end{array}
\right.
\end{equation}
where $(\gamma, w)$ denote the perturbations and $(\hat \gamma, \hat w)= ((\gamma_n)_{n\in\Z}, (w_n)_{n\in\Z})$ are those Fourier modes. 
Now let us suppose that $(\psi, \omega)$ in \eqref{solution} are smooth enough and really solutions of \eqref{VS2} for the moment. Since $\Delta_{r,\theta}\psi_*=-\omega_*$ and $\Delta_{r,\theta}\omega_*=\Delta_{r,\theta}\phi_*$, the perturbations $(\hat \gamma, \hat w)$ should satisfy the following system.
\begin{equation}
\label{VSp}
\left\{
\begin{array}{ll}
   \displaystyle \Delta_{r,n}\gamma_n= -w_n, \hspace{5pt}  &n\in\Z,\ r\geq 0, \\
  \displaystyle  \Delta_{r,n}w_n=
\G^*_n
+\Delta_{r,n} \varphi_n, \hspace{5pt}  &n\in\Z,\ r\geq 0.
  \end{array}
\right.
\end{equation}
Here $\hat \varphi=(\varphi_n)_{n\in\Z}$ denotes the Fourier mode of $\varphi$, $\Delta_{r,n}:=\partial^2_r+(1/r)\dr-(n^2/r^2)I$, and
\begin{equation}
\label{Gn1}
\G^*_n=\G_n(\psi_*, \hat\gamma):=
-\frac{1}{r^2}\dr(r\dr \D_n(\psi_*, \hat\gamma)+ in \E_n(\psi_*, \hat\gamma)) +\frac{1}{r^3}(1-n^2)\D_n(\psi_*, \hat\gamma)
\end{equation}
is the Fourier mode of $G$ associated with \eqref{G1}, where
\begin{alignat}{2}
&\D_n(\psi_*, \hat\gamma):= i\sum_{k+l=n}k\gamma_k\dr \gamma_l
+ in\gamma_n \dr \psi_*,& &n\in\Z, \label{Dn}\\
&\E_n(\psi_*, \hat\gamma):=-\frac1r \sum_{k+l=n}kl\gamma_k\gamma_l
-r\sum_{k+l=n}\dr\gamma_k\dr \gamma_l
-2r\dr\gamma_n \dr \psi_*,&\hspace{10pt} &n\in\Z\backslash\{0\}. \label{En}
\end{alignat}
On the other hand, in the notation of \eqref{G2}, $\G^*_n$ is also expressed as
\begin{equation}
\label{Gn2}
\G^*_n=\H_n(\psi_*, \omega_*, \hat\gamma, \hat w):=\frac{i}{r} \sum_{k+l=n}(k\gamma_k\dr w_l-lw_l\dr \gamma_k)+\frac{in\dr\omega_*}{r}\gamma_n-\frac{in\dr\psi_*}{r}w_n. 
\end{equation}
Since $f(r)=r^{\pm |n|}$ are fundamental solutions of the ordinary equation $\Delta_{r,n}f=0$ for each $n\in\Z\backslash\{0\}$, $(\hat \gamma, \hat w)$ satisfying the system \eqref{VSp} are expressed as
\renewcommand{\arraystretch}{1.5}
\begin{align}
\label{gamma}
&\gamma_n(r)=
\left\{
\begin{array}{ll}
   \displaystyle I^\infty_{|n|}[w_n](r)+ J^0_{|n|}[w_n](r),  &r\geq0,\ n\in\Z\backslash\{0\},\\
   \displaystyle -\int^r_0\frac1s\int_0^s tw_0(t)dtds,  &r\geq0, \ n=0,
  \end{array}
\right.
\\
\label{w}
&w_n(r)=
\left\{
\begin{array}{ll}
  \displaystyle   -I^\infty_{|n|}[\G^*_n](r)- J^0_{|n|}[\G^*_n](r)+\varphi_n(r),  &r\geq0,\ n\in\Z\backslash\{0\},\\
  \displaystyle  -\int^\infty_r\frac1s\int_0^s t\G^*_0(t)dtds+\varphi_0(r),  &r\geq0, \ n=0.
  \end{array}
\right.
\end{align}
Here we define integrations $I^T_z$ and $J^t_z$ as
\begin{align*}
&I^T_z[f](r):= \frac{r^z}{2z}\int^T_r s^{1-z}f(s)ds,\hspace{5pt}z\in\Com\backslash\{0\},\ 0\leq r\leq T\  (0\leq r<\infty\ {\rm if}\ T=\infty),\\
&J^t_z[f](r):= \frac{1}{2zr^z}\int^r_t s^{1+z}f(s)ds,\hspace{5pt}z\in\Com\backslash\{0\},\ t\leq r<\infty.
\end{align*}

\renewcommand{\arraystretch}{1}
We will mainly analyze the expressions \eqref{gamma} and \eqref{w} with $\G^*_n$ of the expression \eqref{Gn1},
while we will apply \eqref{Gn2} to confirm the non-singularity at $r=0$.
It should be emphasized here that the expression \eqref{Gn1} reveals the key cancellation for $|n|=1$ in achieving the desired spatial decay, which is difficult to see if one uses only \eqref{Gn2}.
However, there exists a problem in analysis of \eqref{w} when $|n|= 2$.
Indeed, it is difficult to derive the decay property $\lim_{r\to\infty}(1+r)^2|w_{\pm2}(r)|=0$ from the expression \eqref{w}.
Therefore, in order to solve such a problem for $n=\pm2$, we utilize the effect of $\phi_*$ as follows.
Since
\begin{equation*}
\dr\omega_*(r)= \dr\phi_*(r)\equiv 0,\hspace{10pt}
\dr\psi_*(r)=-\frac1r\int_0^r s\phi_*(s)ds =-\frac{\mu_*}{r}
\end{equation*}
for every $r\geq R_*$, we can rewrite the system \eqref{VSp} locally as 
\begin{equation*}
\left\{
\begin{array}{ll}
    \Delta_{r,n}\gamma_n= -w_n, \hspace{5pt}  &n\in\Z,\ r> R_*,\\
    \Delta_{r,\zeta_n}w_n=\G^0_n
+\Delta_{r,n} \varphi_n, \hspace{5pt}  &n\in\Z,\ r> R_*,
  \end{array}
\right.
\end{equation*}
where $\zeta_n:=(n^2+in\mu_*)^{\frac12}$, $\Delta_{r,\zeta_n}:=\partial^2_r+(1/r)\dr-(\zeta_n^2/r^2)Id$,
and
\begin{equation}
\label{g0n}
\G^0_n
:=\G_n(0, \hat\gamma)=\H_n(0,0,\hat\gamma, \hat w).
\end{equation}
Then using the fundamental solutions $f(r)=r^{\pm \zeta_n}$ for $\Delta_{r,\zeta_n}f=0$, we have the other expression of $w_n$ such that
\begin{equation}
\label{w2}
w_n(r)=-I^\infty_{\zeta_n}[\G^0_n+\Delta_{r,n} \varphi_n](r)- J^0_{\zeta_n}[\G^0_n+\Delta_{r,n} \varphi_n](r),\hspace{5pt}r> R_*,\ n\in\Z\backslash\{0\}.
\end{equation}
Actually, we can find the better decay property for $w_{\pm2}$ of \eqref{w2} than those of \eqref{w}.
Hence in the case of $|n|=2$, we set $w_n$ as
\begin{equation}
\label{w3}
w_n(r)=
\left\{
\begin{array}{ll}
   \displaystyle -I^{R_*}_{|n|}[\G^*_n](r)- J^0_{|n|}[\G^*_n](r)+\varphi_n(r)+c_1r^{|n|}, &0\leq r\leq R_*, \\
  \displaystyle -I^\infty_{\zeta_n}[\G^0_n+\Delta_{r,n} \varphi_n](r)- J^{R_*}_{\zeta_n}[\G^0_n+\Delta_{r,n} \varphi_n](r)+c_2r^{-\zeta_n}, &R_*\leq r<\infty,
  \end{array}
\right.
\end{equation}
where $c_1, c_2\in\R$ are constants.
In this case, we should fix $c_1$ and $c_2$ so that $w_n$ defined by \eqref{w3} becomes continuous and differentiable at $r=R_*$, the detail of which will be stated later.

\subsection{Linear estimates for \eqref{gamma}}

Based on the above discussion, let us construct the theory to obtain solutions of \eqref{VS2} and the original system \eqref{NS}.
 In what follows, we introduce the notation $A \lesssim B$, which means that there exists a constant $c>0$ independent of any given data such that $A\leq cB$.
First of all, let us define the linear map $\L$ associated with \eqref{gamma}, which has the following property.

\begin{prop}\label{prop:1}
Let $0<\alpha<1/2$ and $\kappa>1$. For every $\hat w\in \U^1_{\alpha+2,\kappa+2}$,
define the map $\L : \hat w \mapsto \hat\gamma$ as
\begin{equation*}
\label{a}
\gamma_n(r):=
\left\{
\begin{array}{ll}
   \displaystyle I^\infty_{|n|}[w_n](r)+ J^0_{|n|}[w_n](r), \  &n\in\Z\backslash\{0\},\\
   \displaystyle -\int^r_0\frac1s\int_0^s tw_0(t)dtds, \  &n=0.
  \end{array}
\right.
\end{equation*}
Then $\hat\gamma$ belongs to $\V^{2}_{\alpha,\kappa+4}$, and there hold
\begin{equation}
\label{L}
\|\gamma_0\|_{\V^{2}_0}\ls \frac{1}{\alpha} \|\hat w\|_{{\U}^0_{\alpha+2,\kappa+2}},\hspace{10pt}
\|\tilde \gamma \|_{\tilde{\V}^{2}_{\alpha,\kappa+4}}\ls \|\hat w\|_{{\U}^1_{\alpha+2,\kappa+2}}.
\end{equation}
More precisely, we see
\begin{alignat}{2}
 \label{L2}
&\M^l_{0, 0}[\gamma_0](r)\ls\|\hat w\|_{\U^0_{\alpha+2, \kappa+2}}\min\{\alpha^{-1}, r^{2-l}\}& ,
\hspace{10pt} &l=1,2,\\
\label{L22}
&\M^l_{\alpha, \kappa+4}[\gamma_n](r)\ls
\left\{
\begin{array}{lll}
\displaystyle \|\hat w\|_{\U^0_{\alpha+2, \kappa+2}}\min\{1, r^{1-l}\},\  &|n|=1\\
\displaystyle \|\hat w\|_{\U^1_{\alpha+2, \kappa+2}}\min\{1, r^{2-l}\},\  &|n|=2\\
\displaystyle \|\hat w\|_{\U^0_{\alpha+2, \kappa+2}}\min\{1, r^{2-l}\},\  &|n|\geq 3
\end{array}
\right\},&
\hspace{10pt} &l=0,1,2.
\end{alignat}
\end{prop}

We note here that the concrete estimate \eqref{L2} and \eqref{L22} are used later to analyze the behaviors of some biliner forms in the neighborhood of $r=0$.\\

{\it Proof of Proposition \ref{prop:1}.}\hspace{5pt}
Let $\hat w\in \U^1_{\alpha+2,\kappa+2}$. By the definition, there holds
\begin{equation}
\label{vnr}
|w_n(r)|= \frac{\M^0_{\alpha+2,\kappa+2}[w_n](r)}{(1+|n|)^{\kappa+2}(1+r)^{\alpha+2}},
\hspace{10pt} n\in\N.
\end{equation}
We first consider the case $n=0$. We have the expression of $\dr^l \gamma_0$, $l=1,2$ for every $r\geq0$ as
\renewcommand{\arraystretch}{2}
\begin{equation*}
\label{eta0}
\begin{array}{ll}
\displaystyle \dr \gamma_0(r)=-\frac1r\int_0^r sw_0(s)ds,\hspace{10pt}
 \partial^2_r \gamma_0(r)=\frac1{r^2}\int_0^r sw_0(s)ds-w_0(r).
\end{array}
\end{equation*} 
These expressions, \eqref{vnr}, and an inequality
\begin{equation*}
\int_0^r \frac{s}{(1+s)^{\alpha+2}}ds \ls \min\left\{\frac1\alpha, r^{2} \right\}
\end{equation*}
yield \eqref{L} for $\gamma_0$ and \eqref{L2}.

We next consider the case $n\neq 0$. Each $\dr^l \gamma_n$, $l=1,2$, is calculated as follows.
\begin{equation*}
\begin{array}{lll}
\displaystyle\dr\gamma_n(r)
= \frac{|n|}{r} I^\infty_{|n|}[w_n](r) -\frac{|n|}{r} J^0_{|n|}[w_n](r),\\
\displaystyle\partial^2_r\gamma_n(r)
= \frac{|n|(|n|-1)}{r^2} I^\infty_{|n|}[w_n](r)+\frac{|n|(|n|+1)}{r^2} J^0_{|n|}[w_n](r)-w_n(r).
\end{array}
\end{equation*}
As for $ I^\infty_{|n|}[w_n]$, we have
\begin{align}
\label{In}
| I^\infty_{|n|}[w_n](r)|\ls
\left\{
\begin{array}{lll}
\|\hat w\|_{\U^0_{\alpha+2,\kappa+2}}r(1+r)^{-(\alpha+1)},\ &|n|=1, \\
\|\hat w\|_{\U^1_{\alpha+2,\kappa+2}}(1+|n|)^{-(\kappa+4)}\min\left\{
r^2,\ r^{-\alpha}
\right\},\ &|n|=2,\\
\|\hat w\|_{\U^0_{\alpha+2,\kappa+2}}(1+|n|)^{-(\kappa+4)}\min\left\{
r^2,\ r^{-\alpha}
\right\},\ &|n|\geq3.
\end{array}
\right. 
\end{align}
Indeed, when $|n|\neq2$, we can derive \eqref{In} by the estimate
\begin{equation*}
| I^\infty_{|n|}[w_n](r)|
\leq \frac{\|\hat w\|_{\U^0_{\alpha+2,\kappa+2}}r^{|n|}}{2|n|(1+|n|)^{\kappa+2}}\int_r^\infty \frac{s^{1-|n|}}{(1+s)^{\alpha+2}}ds,
\end{equation*}
which can be also applied for the case $|n|=2$ with $r\geq1$. For the case $|n|=2$ with $0\leq r<1$, we derive \eqref{In} as
\begin{align*}
| I^\infty_{|n|}[w_n](r)|
&\ls r^2\int_r^1 s^{-1}|w_n(s)|ds+\|\hat w\|_{\U^0_{\alpha+2,\kappa+2}}r^2\int_1^\infty \frac{s^{-1}}{(1+s)^{\alpha+2}}ds\\
&\ls \|\hat w\|_{\U^1_{\alpha+2,\kappa+2}}r^2.
\end{align*}
Here we have used the inequality
\begin{equation}
\label{mean}
|s^{-1}w_n(s)|\leq \sup_{s\geq0}|\dr w_n(s)| \ls \|\hat w\|_{\U^1_{\alpha+2,\kappa+2}},\ n\in\Z\bs\{0\}
\end{equation}
for $0\leq s\leq1$ in the case of $|n|=2$, which is derived by the mean-value theorem. On the other hand, for $J^0_{|n|}[w_n]$, there holds
\begin{align}
\label{Jn}
|J^0_{|n|}[w_n](r)|
&\leq \frac{\|\hat w\|_{\U^0_{\alpha+2,\kappa+2}}r^{-|n|}}{2|n|(1+|n|)^{\kappa+2}}\int_0^r \frac{s^{1+|n|}}{(1+s)^{\alpha+2}}ds\\
&\ls
\|\hat w\|_{\U^0_{\alpha+2,\kappa+2}}(1+|n|)^{-(\kappa+4)}\min\left\{r^2,\ r^{-\alpha}\right\}
,\ |n|\geq1. \nonumber
\end{align}
Here we have used $\alpha<1/2$ when $|n|=1$. Using these estimates for integrals and the expressions of $\partial^l_r\gamma_n$ above, we obtain \eqref{L} and \eqref{L22}. \qed\\

\renewcommand{\arraystretch}{1}
\subsection{Bilinear estimates of $\D_n$, $\E_n$, and $\H_n$}

Using Proposition \ref{prop:1} and the inequality
\begin{equation}
\label{sum}
\sum_{k+l=n}\frac{1}{(1+|k|)^{\lambda_1}(1+|l|)^{\lambda_2}}\ls \frac{1}{(1+|n|)^{\min\{\lambda_1,\lambda_2\}}},
\hspace{10pt}n\in\Z,\ \lambda_1, \lambda_2>1
\end{equation}
(see Hillairet-Wittwer \cite[Section 5]{HW2013} for example),
we can easily show the following estimates for the bilinear forms $\D_n$ and $\E_n$, which compose $\G^*_n$, as follows.
\begin{prop}\label{prop:0}
Let $\alpha>0$, $\kappa>1$, and $\hat \gamma\in \V^{2}_{\alpha,\kappa+4}$. Then $\D_n(\psi_*,\hat\gamma)$ and $\E_n(\psi_*,\hat\gamma)$, defined in \eqref{Dn} and \eqref{En}, respectively, are estimated as
\begin{align*}
&\|\D_n(\psi_*,\hat\gamma)\|_{\U^1_{\alpha+1, \kappa+3}} \ls 
\left\{
\begin{array}{lll}
\displaystyle \|\tilde \gamma \|^2_{\tilde{\V}^{2}_{\alpha,\kappa+4}} , &n=0,\\
\displaystyle \|\hat \gamma \|_{{\V}^{2}_{\alpha,\kappa+4}}\|\tilde \gamma \|_{\tilde{\V}^{2}_{\alpha,\kappa+4}}+\|\psi_*\|_{\V^{2}_0}\|\tilde \gamma \|_{\tilde{\V}^{2}_{\alpha,\kappa+4}}&n\neq 0,
\end{array}
\right. \\
&\|\E_n(\psi_*,\hat\gamma)\|_{\U^1_{\alpha+1, \kappa+3}} \ls
\|\hat \gamma \|_{{\V}^{2}_{\alpha,\kappa+4}}\|\tilde \gamma \|_{\tilde{\V}^{2}_{\alpha,\kappa+4}}+\|\psi_*\|_{\V^{2}_0}\|\tilde \gamma \|_{\tilde{\V}^{2}_{\alpha,\kappa+4}},  \hspace{23pt}n\neq 0.
\end{align*}
\end{prop}
Combining Proposition \ref{prop:1}-\ref{prop:0}, we obtain the following corollary.

\begin{cor} \label{cor:1}
Let $0<\alpha<1/2$, $\kappa>1$, $\hat w\in  \U^1_{\alpha+2, \kappa+2}$, and let $\hat\gamma=\L(\hat w)$.

$(1)$The $\U^1_{\alpha+1, \kappa+3}$ norms of $\D_n(\psi_*,\hat\gamma)$ and $\E_n(\psi_*,\hat\gamma)$ are estimated as
\begin{align*}
&\|\D_n(\psi_*,\hat\gamma)\|_{\U^1_{\alpha+1, \kappa+3}} \ls 
\left\{
\begin{array}{lll}
\displaystyle \|\hat w\|^2_{\U^1_{\alpha+2, \kappa+2}}, &n=0,\\
\displaystyle \left(\frac1\alpha\|\hat w\|^2_{\U^1_{\alpha+2, \kappa+2}}+\|\psi_*\|_{\V^{2}_0}\|\hat w\|_{\U^1_{\alpha+2, \kappa+2}} \right) &n\neq0,
\end{array}
\right.\\
&\|\E_n(\psi_*,\hat\gamma)\|_{\U^1_{\alpha+1, \kappa+3}} \ls 
\left(\frac1\alpha\|\hat w\|^2_{\U^1_{\alpha+2, \kappa+2}}+\|\psi_*\|_{\V^{2}_0}\|\hat w\|_{\U^1_{\alpha+2, \kappa+2}} \right),\hspace{23pt} n\neq0.
\end{align*}

$(2)$ There hold
\begin{equation}
\label{dnat0}
\D_n(0)=\E_n(0)=0,\ \ n\in\Z\bs\{0\}.
\end{equation}
For the case $n=0$, it holds for $0\leq r\leq 1$ that
\begin{equation}
\label{dn0at0}
|\partial^l_r \D_0(r)| \ls \|\hat w\|^2_{\U^1_{\alpha+2, \kappa+2}}r^{2-l},
\ \ l=0,1.
\end{equation}
\end{cor}

{\it Proof of Corollary \ref{cor:1}.}\hspace{5pt}
(1) is a direct result from \eqref{L} and Proposition \ref{prop:0}. 
We can also show \eqref{dnat0} easily by using \eqref{L2} and \eqref{L22}, since they may hold $\partial_r\gamma_n(0)\neq 0$ only if $n=\pm1$.
For $\D_0$, which is expressed as
\begin{equation*}
\D_0(r)= i\sum_{k\in\Z\backslash\{0\}}k\gamma_k(r)\dr \gamma_{-k}(r),
\end{equation*}
we can show \eqref{dn0at0} as follows. Since $ \overline{\gamma_k}=\gamma_{-k}$ and hence $\gamma_{-k}\dr \gamma_k= \overline{\gamma_k\dr \gamma_{-k}}$, we see
\begin{equation*}
\D_0=i\sum_{k=1}^\infty k(\gamma_k\dr \gamma_{-k}-\gamma_{-k}\dr \gamma_k)
=-2\sum_{k=1}^\infty k\Im[\gamma_k\dr \gamma_{-k}],
\end{equation*}
where $\Im[f]$ denotes the imaginary part of $f$. Moreover, since
\begin{equation*}
\dr \gamma_{-k}=\frac{|k|}{r} I^\infty_{|k|}[\overline{w_k}] -\frac{|k|}{r} J^0_{|k|}[\overline{w_k}]
=\frac{|k|}{r}( \overline{I^\infty_{|k|}[w_k]} - \overline{J^0_{|k|}[w_k]}),
\end{equation*}
we obtain the detailed expression of $\D_0$ as
\begin{equation*}
\D_0=\frac2r \sum_{k=1}^\infty k^2 \Im[I^\infty_k[w_k]\overline{J^0_k[w_k]}-\overline{I^\infty_k[w_k]}J^0_k[w_k]].
\end{equation*}
Therefore, from the property \eqref{In} and \eqref{Jn}, we have \eqref{dn0at0}.\qed\\

Actually, we will use \eqref{dn0at0} to compute the zero mode of vorticity in the neighborhood of $r=0$. For analysis of the non-zero mode with small $r$, we use $\H_n$ instead of $\G^*_n$. Here we state the estimate of $\H_n$ as follows.

\begin{prop}\label{prop:2}
Let $0<\alpha<1/2$, $\kappa>1$, $\hat w\in  \U^1_{\alpha+2, \kappa+2}$, and let $\hat\gamma=\L(\hat w)$. Then $\H_n(\psi_*, \omega_*, \hat\gamma, \hat w)$ defined in \eqref{Gn2} is estimated as
\begin{equation}
\label{H}
|\H_n(r)|\ls \frac1{(1+|n|)^{\kappa+1}}\left\{ \frac{\alpha^{-1}\|\hat w\|^2_{\U^1_{\alpha+2,\kappa+2}}
+\nu_*\|\hat w\|_{\U^1_{\alpha+2,\kappa+2}}}{(1+r)^{\alpha+2}} 
+\frac{\|\hat w\|^2_{\U^1_{\alpha+2,\kappa+2}}}{(1+r)^{\alpha+1}}
\right\}
\end{equation}
for every $0\leq r<\infty$ and $n\in\Z$.
\end{prop}

\vspace{5pt}

{\it Proof of Proposition \ref{prop:2}.}\hspace{5pt}
We decompose $\H_n$ as
\begin{align*}
\H_n
&=\frac{i}{r} \sum_{k+l=n}(k\gamma_k\dr w_l)-\frac{i}{r} \sum_{k+l=n}(lw_l\dr \gamma_k)+\frac{in\dr\omega_*}{r}\gamma_n-\frac{in\dr\psi_*}{r}w_n\\
&=:\H_{n,1}-\H_{n,2}+\H_{n,3}-\H_{n,4}.
\end{align*}
By \eqref{L22}, the inequality \eqref{sum} and 
\begin{equation*}
|\dr w_n(r)|= \frac{\M^1_{\alpha+2,\kappa+2}[w_n](r)}{(1+|n|)^{\kappa+1}(1+r)^{\alpha+3}},
\hspace{10pt} n\in\Z,
\end{equation*}
we easily see that
\begin{equation*}
|\H_{n,1}(r)|\ls 
\frac{\|\hat w\|^2_{\U^1_{\alpha+2,\kappa+2}}}{(1+|n|)^{\kappa+1}(1+r)^{\alpha+3}}.
\end{equation*}
Similarly, we have
\begin{equation*}
|\H_{n,3}(r)|\ls
\frac{\|\hat w\|_{\U^1_{\alpha+2,\kappa+2}}(1+r)^2|\dr\omega_*(r)|}{(1+|n|)^{\kappa+1}(1+r)^{\alpha+2}}
\ls \frac{\|\hat w\|_{\U^1_{\alpha+2,\kappa+2}}\nu_*}{(1+|n|)^{\kappa+1}(1+r)^{\alpha+2}}.
\end{equation*}
On the other hand, if $|k|\neq 1$, then by \eqref{L2} and \eqref{L22}, we see that
\begin{equation*}
\left| lw_l(r)\frac{\dr\gamma_k(r)}{r}\right|
\ls \frac{\alpha^{-1}\|\hat w\|^2_{\U^1_{\alpha+2,\kappa+2}}}{(1+|l|)^{\kappa+1}(1+|k|)^{\kappa+3}(1+r)^{\alpha+2}},
\end{equation*}
while if $|k|=1$, then by \eqref{mean}, we see
\begin{equation*}
\left| l\frac{w_l(r)}{r}\dr\gamma_k(r)\right|
\ls \frac{\|\hat w\|^2_{\U^1_{\alpha+2,\kappa+2}}}{(1+|l|)^{\kappa+1}(1+|k|)^{\kappa+3}(1+r)^{\alpha+1}}.
\end{equation*}
Therefore, using \eqref{sum}, we have
\begin{equation*}
|\H_{n,2}(r)|\ls 
\frac{\|\hat w\|^2_{\U^1_{\alpha+2,\kappa+2}}}{(1+|n|)^{\kappa+1}}
\left(
\frac1{(1+r)^{\alpha+2}}+\frac1{(1+r)^{\alpha+1}}
\right).
\end{equation*}
As for $\H_{n,4}$, we see
\begin{equation*}
|\H_{n,4}(r)|\ls 
\frac{\|\hat w\|_{\U^1_{\alpha+2,\kappa+2}}(r^{-1}\dr\psi_*(r))}{(1+|n|)^{\kappa+1}(1+r)^{\alpha+3}}
\ls \frac{\|\hat w\|_{\U^1_{\alpha+2,\kappa+2}}\nu_*}{(1+|n|)^{\kappa+1}(1+r)^{\alpha+3}}.
\end{equation*}
Here we have used the inequality
\begin{equation*}
\left|\frac1{r^2}\int_0^r s\phi_*(s)ds\right| \ls \sup_{0\leq r\leq R_*}|\phi_*(r)|<\nu_*.
\end{equation*}
Together with the above estimates, we obtain \eqref{H}. \qed\\

\renewcommand{\arraystretch}{1.5}
\subsection{Solution map for \eqref{w} and \eqref{w3}}
Finally, let us define the map $\S$ associated with \eqref{w} and \eqref{w3} as follows. 
For fixed $0<\alpha<1/2$ and $\kappa>1$, we take $(\hat w, \hat\sigma)\in \U^1_{\alpha+2,\kappa+2} \times \U^1_{\alpha+2,\kappa+2}$ arbitrarily, and let $\hat\gamma=\L(\hat w)$.
Then we define 
$\hat y=\S(\hat w, \hat\sigma)$ as
\begin{equation}
\label{defyn}
y_n(r):=
\left\{
\begin{array}{llll}
 \displaystyle    -I^\infty_{|n|}[\G^*_n](r)- J^0_{|n|}[\G^*_n](r)+\sigma_n(r),  &r\geq0,\ n\in\Z\backslash\{0,\pm2\},\\
\displaystyle  y_{n,1}(r):=\P_{n,1}(r)+ \Q_{n,1}[\sigma_n](r)
&0\leq r\leq R_*,\ n=\pm2,\\
\displaystyle y_{n,2}(r):=\P_{n,2}(r)+ \Q_{n,2}[\sigma_n](r)
&r> R_*,\ n=\pm2,\\
\displaystyle   -\int^\infty_r\frac1s\int_0^s t\G^*_0(t)dtds+\sigma_n(r),  &r\geq0, \ n=0,
  \end{array}
\right.
\end{equation}
where each $\G^*_n$ is that of \eqref{Gn1} or \eqref{Gn2},
\begin{align*}
\P_{n,1}(r)
:=&-I^{R_*}_{|n|}[\G^*_n](r)- J^0_{|n|}[\G^*_n](r)\\
&\hspace{5pt} +\frac{(r/R_*)^{|n|}}{|n|+\zeta_n}
\left\{(\zeta_n-|n|)J^0_{|n|}[\G^*_n](R_*)-2\zeta_n I^\infty_{\zeta_n}[\G^0_n](R_*)\right\},\ 0\leq r\leq R_*,\\
\P_{n,2}(r)
:=&-I^\infty_{\zeta_n}[\G^0_n](r)- J^{R_*}_{\zeta_n}[\G^0_n](r) \\
&\hspace{5pt}+\frac{(R_*/r)^{\zeta_n}}{|n|+\zeta_n}
\left\{-(\zeta_n-|n|)I^\infty_{\zeta_n}[\G^0_n](R_*)-2|n|J^0_{|n|}[\G^*_n](R_*)\right\},\ r>R_*
\end{align*}
with $\G^0_n$ of \eqref{g0n}, and
\begin{align*}
\Q_{n,1}[\sigma_n](r)
:=&\sigma_n(r) \\
&\hspace{5pt} +\frac{(r/R_*)^{|n|}}{|n|+\zeta_n}
\left\{-2\zeta_n I^\infty_{\zeta_n}[\Delta_{r,n}\sigma_n](R_*)-\zeta_n\sigma_n(R_*)-R_*\dr\sigma_n(R_*)\right\},\ 0\leq r\leq R_*,\\
\Q_{n,2}[\sigma_n](r)
:=&-I^\infty_{\zeta_n}[\Delta_{r,n}\sigma_n](r)- J^{R_*}_{\zeta_n}[\Delta_{r,n}\sigma_n](r) \\
&\hspace{5pt}+\frac{(R_*/r)^{\zeta_n}}{|n|+\zeta_n}
\left\{-(\zeta_n-|n|)I^\infty_{\zeta_n}[\Delta_{r,n}\sigma_n](R_*)+|n|\sigma_n(R_*)-R_*\dr\sigma_n(R_*)\right\},
\ r>R_*.\
\end{align*}
In the case $|n|=2$, we see from the above expressions that if $y_{n,1}$ and $y_{n,2}$ are well-defined and differentiable in each domain, there hold
\begin{equation}
\label{match1}
\partial^l_r \P_{n,1}(R_*)=\lim_{r\to R_*+0}\partial^l_r\P_{n,2}(r),
\hspace{10pt}l=0,1
\end{equation}
and
\begin{equation}
\label{match2}
\partial^l_r\Q_{n,1}[\sigma_n](R_*)=\lim_{r\to R_*+0}\partial^l_r\Q_{n,2}[\sigma_n](r),\hspace{10pt}l=0,1,
\end{equation}
so that each of $y_{\pm2}$ belongs to the $C^1$ class on $[0,\infty)$.

First of all, we compute the integrals associated with $\G^*_n$ and $\G^0_{\pm2}$. 

\begin{prop}
\label{prop:3}
Suppose that $\kappa>1$ and $\hat w \in \U^1_{\alpha+2, \kappa+2}$.

$(1)$ Let $n\in \Z\bs\{0, \pm2\}$ and $0<\alpha<1/2$. 
Then for every $0\leq r<\infty$, there holds
\begin{equation*}
\M^0_{\alpha+2,\kappa+2}[I^{\infty}_{|n|}[\G^*_n]](r)
+\M^0_{\alpha+2,\kappa+2}[J^0_{|n|}[\G^*_n]](r)
\ls \frac1\alpha\|\hat w\|^2_{\U^1_{\alpha+2, \kappa+2}}+\nu_*\|\hat w\|_{\U^1_{\alpha+2, \kappa+2}}.
\end{equation*}

$(2)$ Let $n=\pm2$, $0<\alpha<1/2$ and $0\leq r\leq R_*$. Then 
\begin{equation*}
\M^0_{\alpha+2,\kappa+2}[I^{R_*}_{|n|}[\G^*_n]](r)
\ls \frac1\alpha\|\hat w\|^2_{\U^1_{\alpha+2, \kappa+2}}+\nu_*\|\hat w\|_{\U^0_{\alpha+2, \kappa+2}},
\end{equation*}
\begin{equation*}
\M^0_{\alpha+2,\kappa+2}[J^0_{|n|}[\G^*_n]](r)
\ls R_*^\alpha
\left(
\frac1\alpha\|\hat w\|^2_{\U^1_{\alpha+2, \kappa+2}}+\nu_*\|\hat w\|_{\U^1_{\alpha+2, \kappa+2}}\right).
\end{equation*}

$(3)$ Let $n=\pm2$, $0<\alpha<\min\{1/2, \rho_*\}$ and $R_*\leq r<\infty$. Then
\begin{equation*}
 \M^0_{\alpha+2,\kappa+2}[I^\infty_{\zeta_n}[\G^0_n]](r)
\ls \frac1\alpha\|\hat w\|^2_{\U^1_{\alpha+2, \kappa+2}},
\end{equation*}
\begin{equation*}
\M^0_{\alpha+2,\kappa+2}[J^{R_*}_{\zeta_n}[\G^0_n]](r)
\ls \frac1{\alpha(\rho_*-\alpha)}\|\hat w\|^2_{\U^1_{\alpha+2, \kappa+2}}.
\end{equation*}
\end{prop}

{\it Proof of Proposition \ref{prop:3}.}\hspace{5pt}
(1) We firstly consider the case $1\leq r<\infty$. Using \eqref{Gn1} and integration by parts, we see
\begin{align*}
I^\infty_{|n|}[\G^*_n](r)
&=-\frac{r^{|n|}}{2|n|}\int_r^\infty \left[ s^{-1-|n|}\ds(s\ds\D_n(s)+in\E_n(s)) + s^{-2-|n|}(1-n^2) \D_n(s)\right]ds\\
&=\frac1{2|n|}\dr\D_n(r)+\frac{1+|n|}{2|n|r}\D_n(r)+\frac{in}{2|n|r}\E_n(r) \\
&\hspace{10pt} -(1+|n|)r^{|n|}\int_r^\infty s^{-2-|n|}\D_n(s)ds -\frac{in(1+|n|)}{2|n|}r^{|n|}\int_r^\infty s^{-2-|n|}\E_n(s)ds
\end{align*}
and
\begin{align*}
J^0_{|n|}[\G^*_n](r)
&=-\frac{1}{2|n|r^{|n|}}\int_0^r \left[ s^{|n|-1}\ds(s\ds\D_n(s)+in\E_n(s)) + s^{|n|-2}(1-n^2) \D_n(s)\right]ds\\
&=-\frac1{2|n|}\dr\D_n(r)+\frac{|n|-1}{2|n|r}\D_n(r)-\frac{in}{2|n|r}\E_n(r) \\
&\hspace{10pt} +\frac{1-|n|}{r^{|n|}}\int_0^r s^{|n|-2}\D_n(s)ds +\frac{in(|n|-1)}{2|n|r^{|n|}}\int_0^r s^{|n|-2}\E_n(s)ds,
\end{align*}
where $\D_n=\D_n(\psi_*, \hat\gamma)$ and $\E_n=\E_n(\psi_*, \hat\gamma)$, and we have used \eqref{dnat0} when $|n|=1$.
Here we should note that if $|n|=1$, then the second, fourth and fifth parts of $J^0_{|n|}[\G^*_n]$ are vanished.
Corollary \ref{cor:1} (1) yield
\begin{align*}
r^{|n|}\int_r^\infty s^{-2-|n|}|\D_n(s)|ds
&=(1+|n|)^{-(\kappa+3)}r^{|n|}\int_r^\infty \frac{\M^0_{\alpha+1, \kappa+3}[\D_n](s)}{s^{2+|n|}(1+s)^{\alpha+1}}ds\\
&\ls (1+|n|)^{-(\kappa+4)}\left(\frac1\alpha\|\hat w\|^2_{\U^1_{\alpha+2, \kappa+2}}+\|\psi_*\|_{\V^{2}_0}\|\hat w\|_{\U^1_{\alpha+2, \kappa+2}} \right)r^{-(\alpha+2)}
\end{align*}
and the same estimate for $\E_n$ for every $n\in \Z\bs\{0\}$. 
Hence, we obtain
\begin{equation}
\label{Iinfty1}
\M^0_{\alpha+2,\kappa+2}[I^\infty_{|n|}[\G^*_n]](r)
\ls \frac1\alpha\|\hat w\|^2_{\U^1_{\alpha+2, \kappa+2}}+\|\psi_*\|_{\V^{2}_0}\|\hat w\|_{\U^1_{\alpha+2, \kappa+2}},
\ \ n\in\Z\bs\{0\}.
\end{equation}
for every $1\leq r<\infty$. 
On the other hand, it holds that
\begin{align}
\label{jn0}
&\frac1{r^{|n|}}\int_0^r s^{|n|-2}|\D_n(s)|ds
=\frac{(1+|n|)^{-(\kappa+3)}}{r^{|n|}}\int_0^r \frac{s^{|n|-2}}{(1+s)^{\alpha+1}}\M^0_{\alpha+1, \kappa+3}[\D_n](s)ds \\
&\hspace{10pt}\ls (1+|n|)^{-(\kappa+4)}\left(\frac1\alpha\|\hat w\|^2_{\U^1_{\alpha+2, \kappa+2}}+\|\psi_*\|_{\V^{1}_0}\|\hat w\|_{\U^1_{\alpha+2, \kappa+2}} \right)
\times
\left\{
\begin{array}{ll}
r^{-(\alpha+2)}, |n|\geq3,\\
r^{-2}, |n|=2,
\end{array}
\right. \nonumber
\end{align}
where we should use $0<\alpha<1/2$ when $|n|=3$. 
Since we also obtain the same estimate for $\E_n$, we have
\begin{equation}
\label{Iinfty2}
\M^0_{\alpha+2,\kappa+2}[J^0_{|n|}[\G^*_n]](r)
\ls \frac1\alpha\|\hat w\|^2_{\U^1_{\alpha+2, \kappa+2}}+\|\psi_*\|_{\V^{2}_0}\|\hat w\|_{\U^1_{\alpha+2, \kappa+2}},\hspace{5pt}  n\in\Z\bs\{0,\pm2\}
\end{equation}
for every $1\leq r<\infty$. 

We next consider the case $0\leq r\leq 1$. 
Using \eqref{Gn2} and \eqref{H}, we have
\begin{align}
\label{HH1}
|I^{\infty}_{|n|}[\G^*_n](r)|
&\leq 
\frac{r^{|n|}}{2|n|}\int_r^\infty s^{1-|n|}|\H_n(s)|ds \\
&\ls
\frac1{(1+|n|)^{\kappa+2}}\left(
\frac1\alpha \|\hat w\|^2_{\U^1_{\alpha+2,\kappa+2}}
+\nu_*\|\hat w\|_{\U^1_{\alpha+2,\kappa+2}}
\right)r\nonumber
\end{align}
and
\begin{align}
\label{HH2}
|J^{0}_{|n|}[\G^*_n](r)|
&\leq 
\frac{1}{2|n|r^{|n|}}\int_0^r s^{1+|n|}|\H_n(s)|ds \\
&\ls
\frac1{(1+|n|)^{\kappa+2}}\left(
\frac1\alpha \|\hat w\|^2_{\U^1_{\alpha+2,\kappa+2}}
+\nu_*\|\hat w\|_{\U^1_{\alpha+2,\kappa+2}}
\right)r^2,\nonumber
\end{align}
which, \eqref{Iinfty1}, \eqref{Iinfty2}, and the inequality $\|\psi_*\|_{\V^{2}_0}\ls\nu_*$ yield (1).

(2) For $I^{R_*}_{|n|}[\G^*_n]$, which is expanded as
\begin{align*}
I^{R_*}_{|n|}[\G^*_n](r)
&=\frac1{2|n|}\dr\D_n(r)+\frac{1+|n|}{2|n|r}\D_n(r)+\frac{in}{2|n|r}\E_n(r) \\
&\hspace{10pt} -(1+|n|)r^{|n|}\int_r^{R_*} s^{-2-|n|}\D_n(s)ds -\frac{in(1+|n|)}{2|n|}r^{|n|}\int_r^{R_*} s^{-2-|n|}\E_n(s)ds\\
&\hspace{10pt} -\frac{r^{|n|}}{2|n|}\left(
R_*^{-|n|}\dr\D_n(R_*) +(1+|n|)\D_n(R_*) +inR_*^{-|n|-1}\E_n(R_*)
\right),\ 1\leq r\leq R_*,
\end{align*}
we easily obtain the same estimate as \eqref{Iinfty1} for $1\leq r\leq R_*$, that is,
\begin{equation*}
\M^0_{\alpha+2,\kappa+2}[I^{R_*}_{|n|}[\G^*_n]](r)
\ls \frac1\alpha\|\hat w\|^2_{\U^1_{\alpha+2, \kappa+2}}+\|\psi_*\|_{\V^{2}_0}\|\hat w\|_{\U^1_{\alpha+2, \kappa+2}},\hspace{5pt}1\leq r\leq R_*.
\end{equation*}
Moreover, by using \eqref{jn0} with $|n|=2$ and the inequality $(1+r)^{-2}\ls R_*^\alpha(1+r)^{-2-\alpha}$ in $1\leq r\leq R_*$, we also have the estimate for $J^0_{|n|}[\G^*_n]$ as
\begin{equation*}
\M^0_{\alpha+2,\kappa+2}[J^0_{|n|}[\G^*_n]](r)
\ls R_*^\alpha
\left(
\frac1\alpha\|\hat w\|^2_{\U^1_{\alpha+2, \kappa+2}}+\|\psi_*\|_{\V^{2}_0}\|\hat w\|_{\U^1_{\alpha+2, \kappa+2}}\right) ,\hspace{5pt}1\leq r\leq R_*.
\end{equation*}
Moreover, we can also apply the same theory as (1) for the case $0\leq r\leq 1$. This completes the proof of (2).

(3) Let $n=\pm2$ and $r\geq R_*$. First of all, let us confirm that
\begin{align*}
|\zeta_{\pm2}|=2\left[1+\left(\frac{\mu_*}{2}\right)^2\right]^{\frac14},\hspace{10pt}
|r^{\zeta_{\pm2}}|=r^{\Re[\zeta_{\pm2}]}=r^{\rho_*+2},
\end{align*}
where $\Re[f]$ denotes the real part of $f$.

The results of the calculation of $I^\infty_{\zeta_n}[\G^0_n]$ and $J^{R_*}_{\zeta_n}[\G^0_n]$ by integration by parts are as follows.
\begin{align*}
I^\infty_{\zeta_n}[\G^0_n](r)
&=\frac1{2\zeta_n}\dr\D^0_n(r)+\frac{1+\zeta_n}{2\zeta_nr}\D^0_n(r)+\frac{in}{2\zeta_nr}\E^0_n(r) \\
&\hspace{10pt} -\frac{1-n^2-(1+\zeta_n)^2}{2\zeta_n}r^{\zeta_n}\int_r^\infty s^{-2-\zeta_n}\D^0_n(s)ds -\frac{in(1+\zeta_n)}{2\zeta_n}r^{\zeta_n}\int_r^\infty s^{-2-\zeta_n}\E^0_n(s)ds,
\end{align*}
\begin{align*}
J^{R_*}_{\zeta_n}[\G^0_n](r)
&=-\frac1{2\zeta_n}\dr\D^0_n(r)+\frac{\zeta_n-1}{2\zeta_nr}\D^0_n(r)-\frac{in}{2\zeta_nr}\E^0_n(r) \\
&\hspace{10pt} +\frac{1-n^2-(1-\zeta_n)^2}{2\zeta_n r^{\zeta_n}}\int_{R_*}^r s^{\zeta_n-2}\D^0_n(s)ds +\frac{in(\zeta_n-1)}{2\zeta_nr^{\zeta_n}}\int_{R_*}^r s^{\zeta_n-2}\E^0_n(s)ds \\
&\hspace{10pt} +\frac1{2\zeta_n r^{\zeta_n}}\left(
R_*^{\zeta_n}\dr\D^0_n(R_*)-(\zeta_n-1)R_*^{\zeta_n-1}\D^0_n(R_*)+inR_*^{\zeta_n-1}\E^0_n(R_*)
\right).
\end{align*}
Here $\D^0_n:=\D_n(0,\hat\gamma)$ and $\E^0_n:=\E_n(0,\hat\gamma)$ are independent of $\psi_*$. For the integrals, we see
\begin{align*}
\left|r^{\zeta_n}\int_r^\infty s^{-2-\zeta_n}\D^0_n(s)ds\right|
&\ls (1+|n|)^{-(\kappa+3)}r^{\rho_*+2}\int_r^\infty \frac{\M^0_{\alpha+1, \kappa+3}[\D^0_n](s)}{s^{\rho_*+4}(1+s)^{\alpha+1}}ds\\
&\ls (1+|n|)^{-(\kappa+3)}\left(\frac1\alpha\|\hat w\|^2_{\U^1_{\alpha+2, \kappa+2}}\right)r^{-(\alpha+2)},
\end{align*}
\begin{align*}
\left|\frac1{r^{\zeta_n}}\int_{R_*}^r s^{\zeta_n-2}\D^0_n(s)ds\right|
&\ls \frac{(1+|n|)^{-(\kappa+3)}}{r^{\rho_*+2}}\int_{0}^r \frac{s^{\rho_*}}{(1+s)^{\alpha+1}}\M^0_{\alpha+1, \kappa+3}[\D^0_n](s)ds \\
&\ls \frac{(1+|n|)^{-(\kappa+3)}}{\rho_*-\alpha}\left(\frac1\alpha\|\hat w\|^2_{\U^1_{\alpha+2, \kappa+2}} \right)r^{-(\alpha+2)},
\end{align*}
and the same estimates for $\E^0_n$.
Here we have used $0<\alpha<\rho_*$ for the last integration.
Therefore, we obtain the estimate for $I^\infty_{\zeta_n}[\G^0_n]$ and $J^{R_*}_{\zeta_n}[\G^0_n]$ for $n=\pm2$ as (3). \qed\\

Using Proposition \ref{prop:3}, we obtain the following key estimate for the solution map $\S$.
\begin{prop}
\label{prop:4}
Let $0<\alpha<\min\{1/2, \rho_*\}$ and $\kappa>1$. Then the map $\S : (\hat w,  \hat \sigma) \mapsto \hat y$ defined by \eqref{defyn} is bounded from $\U^1_{\alpha+2,\kappa+2} \times  \U^1_{\alpha+2,\kappa+2}$ to $\U^1_{\alpha+2,\kappa+2}$. Moreover, there holds
\begin{equation*}
\label{S}
\|\S(\hat w, \hat\sigma)\|_{\U^1_{\alpha+2,\kappa+2}} \ls
R_*^{\rho_*} \left(\frac1{\alpha(\rho_*-\alpha)} \|\hat w\|^2_{\U^1_{\alpha+2, \kappa+2}}+\nu_*\|\hat w\|_{\U^1_{\alpha+2, \kappa+2}} +\frac{1}{\rho_*-\alpha}\|\hat \sigma\|_{\U^1_{\alpha+2,\kappa+2}}\right).
\end{equation*}
\end{prop}
\vspace{10pt}

{\it Proof of Proposition \ref{prop:4}.}

{\it Step 1: The case $n=0$.} \hspace{5pt}
We see
\begin{equation*}
\G^*_0
=-\frac1{r^2}\dr(r\dr \D_0)+\frac1{r^3}\D_0,
\end{equation*}
where $\D_0=\D_0(0,\hat\gamma)$ depends only on $\hat\gamma$. Using \eqref{dn0at0} and integration by parts, we have
\begin{align*}
y_0(r)
&=-\int_r^\infty \frac1s \int_0^s \left[-\frac1t\dt(t\dt\D_0(t))+\frac1{t^2}\D_0(t)\right]dt +\sigma_0(r)\\
&=-\frac1r\D_0(r)+2 \int_r^\infty \frac{1}{s^2}\D_0(s)ds +\sigma_0(r).\nonumber
\end{align*}
Hence for every $0\leq r<\infty$, we see from \eqref{dn0at0} again that 
\begin{align}
\label{y0}
\M^0_{\alpha+2, 0}[y_0](r)
&\ls
\frac{1+r}{r} \M^0_{\alpha+1, \kappa+3}[\D_0](r)\\
&\hspace{5pt}
+(1+r)^{-(\alpha+2)}\int_r^\infty\frac{\M^0_{\alpha+1, \kappa+3}[\D_0](s)}{s^2(1+s)^{\alpha+1}}ds
+\M^0_{\alpha+2, \kappa+2}[\sigma_0](r) \nonumber\\
&\ls
\|\hat w\|^2_{\U^1_{\alpha+2,\kappa+2}}+\|\hat \sigma\|_{\U^0_{\alpha+2,\kappa+2}}.\nonumber
\end{align}
Moreover, since
\begin{equation*}
\dr y_0(r)
=-\frac1r\dr\D_0(r)-\frac{1}{r^2}\D_0(r)+\dr\sigma_0(r),
\end{equation*}
we see by \eqref{dn0at0} that
\begin{equation}
\label{dry0}
\M^1_{\alpha+2, 0}[y_0](r)
\ls \|\hat w\|^2_{\U^1_{\alpha+2,\kappa+2}}+\|\hat \sigma\|_{\U^1_{\alpha+2,\kappa+2}},\hspace{10pt}
0\leq r<\infty.
\end{equation}
\vspace{5pt}

{\it Step 2: The case $n\in\Z\bs\{0, \pm2\}$.} \hspace{5pt}
By Proposition \ref{prop:3} (1), we obtain
\begin{equation}
\label{ynot2}
\M^0_{\alpha+2, \kappa+2}[y_n](r)
\ls 
\frac1\alpha\|\hat w\|^2_{\U^1_{\alpha+2, \kappa+2}}+\nu_*\|\hat w\|_{\U^1_{\alpha+2, \kappa+2}}+\|\hat \sigma\|_{\U^0_{\alpha+2,\kappa+2}}.
\end{equation}
In addition, since
\begin{equation*}
\dr y_n(r)=-\frac{|n|}{r}I^\infty_{|n|}[\G^*_n](r)+\frac{|n|}{r}J^0_{|n|}[\G^*_n](r)+\dr\sigma_n(r),
\end{equation*}
we also obtain
\begin{equation}
\label{ynot22}
\M^1_{\alpha+2, \kappa+2}[y_n](r)
\ls 
\frac1\alpha\|\hat w\|^2_{\U^1_{\alpha+2, \kappa+2}}+\nu_*\|\hat w\|_{\U^1_{\alpha+2, \kappa+2}}+\|\hat \sigma\|_{\U^1_{\alpha+2,\kappa+2}}.
\end{equation}
(For the case $0\leq r\leq 1$, we pay attention to \eqref{HH1} and \eqref{HH2})
\\

{\it Step 3-1: The case $n=\pm2$ with $1\leq r<\infty$.}\hspace{5pt}
First, let us consider $\P_{n,1}$ and $\P_{n,2}$. 
If $1\leq r\leq R_*$, then
\begin{equation}
\label{fin1}
\begin{array}{ll}
&\displaystyle |I^\infty_{\zeta_n}[\G^0_n](R_*)|
\ls \frac{\M^0_{\alpha+2,\kappa+2}[I^\infty_{\zeta_n}[\G^0_n]](R_*)}{(1+r)^{\alpha+2}(1+|n|)^{\kappa+2}},\\
&\displaystyle |J^0_{|n|}[\G^*_n](R_*)|
\ls \frac{\M^0_{\alpha+2,\kappa+2}[J^0_{|n|}[\G^*_n]](R_*)}{(1+r)^{\alpha+2}(1+|n|)^{\kappa+2}},
\end{array}
\end{equation}
and if $R_*<r<\infty$, then
\begin{equation}
\label{fin2}
\begin{array}{ll}
&\displaystyle \left(\frac{R_*}{r}\right)^{\zeta_n}|I^\infty_{\zeta_n}[\G^0_n](R_*)|
\ls R_*^{\rho_*-\alpha}\frac{\M^0_{\alpha+2,\kappa+2}[I^\infty_{\zeta_n}[\G^0_n]](R_*)}{(1+r)^{\alpha+2}(1+|n|)^{\kappa+2}},\\
&\displaystyle \left(\frac{R_*}{r}\right)^{\zeta_n}|J^0_{|n|}[\G^*_n](R_*)|
\ls R_*^{\rho_*-\alpha}\frac{\M^0_{\alpha+2,\kappa+2}[J^0_{|n|}[\G^*_n]](R_*)}{(1+r)^{\alpha+2}(1+|n|)^{\kappa+2}}.
\end{array}
\end{equation}
Hence, using these and Proposition \ref{prop:3} (2)-(3), we obtain the estimates of $\P_{n,1}$ and $\P_{n,2}$ as
\begin{equation}
\label{P1P2}
\begin{array}{ll}
\displaystyle&\M^0_{\alpha+2,\kappa+2}[\P_{n,1}](r)
\ls R_*^\alpha\left(\frac1\alpha \|\hat w\|^2_{\U^1_{\alpha+2, \kappa+2}}+\nu_*|\hat w\|_{\U^1_{\alpha+2, \kappa+2}} \right),\\
\displaystyle&\M^0_{\alpha+2,\kappa+2}[\P_{n,2}](r)
\ls R_*^{\rho_*} \left(\frac1{\alpha(\rho_*-\alpha)} \|\hat w\|^2_{\U^1_{\alpha+2, \kappa+2}}+\nu_*\|\hat w\|_{\U^1_{\alpha+2, \kappa+2}} \right)
\end{array}
\end{equation}
in $1\leq r\leq R_*$ and $R_*<r<\infty$, respectively.
\renewcommand{\arraystretch}{1}

Secondly, let us consider $\Q_{n,1}$ and $\Q_{n,2}$.
Using integration by parts, we have
\begin{align*}
I^\infty_{\zeta_n}[\Delta_{r,n}\sigma_n](r)
&=\frac{r^{\zeta_n}}{2\zeta_n}\int_r^\infty s^{1-\zeta_n}\left(\frac1s \ds(s\ds\sigma_n(s))-\frac{n^2}{s^2}\sigma_n(s)\right)ds\\
&=-\frac{1}{2\zeta_n}r\dr\sigma_n(r)-\frac12\sigma_n(r)+\left(\frac{\zeta_n}2-\frac2{\zeta_n}\right)r^{\zeta_n}\int_r^\infty s^{-\zeta_n-1}\sigma_n(s)ds
\end{align*}
and
\begin{align*}
J^{R_*}_{\zeta_n}[\Delta_{r,n}\sigma_n](r)
&=\frac1{2\zeta_n r^{\zeta_n}}\int_{R_*}^r s^{1+\zeta_n}\left(\frac1s \ds(s\ds\sigma_n(s))-\frac{n^2}{s^2}\sigma_n(s)\right)ds\\
&=\frac{1}{2\zeta_n}r\dr\sigma_n(r)-\frac12\sigma_n(r)+\left(\frac{\zeta_n}2-\frac2{\zeta_n}\right)\frac{1}{r^{\zeta_n}}\int_{R_*}^r s^{\zeta_n-1}\sigma_n(s)ds\\
&\hspace{10pt}-\frac{R_*^{\zeta_n+1}}{2\zeta_nr^{\zeta_n}}\dr\sigma_n(R_*)
+\frac{R_*^{\zeta_n}}{2r^{\zeta_n}}\sigma_n(R_*).
\end{align*}
Since
\begin{align*}
\left| r^{\zeta_n}\int_r^\infty s^{-\zeta_n-1}\sigma_n(s)ds\right|
&\ls (1+|n|)^{-(\kappa+2)}r^{\rho_*+2}\int_r^\infty\frac{\M^0_{\alpha+2,\kappa+2}[\sigma_n](s)}{s^{\rho_*+3}(1+s)^{\alpha+2}}ds\\
&\ls (1+|n|)^{-(\kappa+2)}\|\hat\sigma\|_{\U^0_{\alpha+2,\kappa+2}}r^{-(\alpha+2)}
\end{align*}
and
\begin{align*}
\left| \frac{1}{r^{\zeta_n}}\int_{R_*}^r s^{\zeta_n-1}\sigma_n(s)ds\right|
&\ls (1+|n|)^{-(\kappa+2)}\frac1{r^{\rho_*+2}}\int_{0}^r\frac{s^{\rho_*+1}\M^0_{\alpha+2,\kappa+2}[\sigma_n](s)}{(1+s)^{\alpha+2}}ds\\
&\ls \frac{(1+|n|)^{-(\kappa+2)}}{\rho_*-\alpha}\|\hat\sigma\|_{\U^0_{\alpha+2,\kappa+2}}r^{-(\alpha+2)},
\end{align*}
we see for $r\geq R_*$ that
\renewcommand{\arraystretch}{1}
\begin{equation}
\label{sigman}
\begin{array}{ll}
\displaystyle \M^0_{\alpha+2,\kappa+2}[I^\infty_{\zeta_n}[\Delta_{r,n}\sigma_n]](r)
\ls \|\hat \sigma\|_{\U^1_{\alpha+2,\kappa+2}},\\
\displaystyle \M^0_{\alpha+2,\kappa+2}[J^{R_*}_{\zeta_n}[\Delta_{r,n}\sigma_n]](r)
\ls \frac{1}{\rho_*-\alpha}\|\hat\sigma\|_{\U^1_{\alpha+2,\kappa+2}}.
\end{array}
\end{equation}
Hence by \eqref{sigman} and a similar discussion to \eqref{fin1} and \eqref{fin2}, $\Q_{n,1}$ and $\Q_{n, 2}$ are estimated as
\begin{equation}
\label{Q1Q2}
\begin{array}{ll}
\displaystyle \M^0_{\alpha+2,\kappa+2}[\Q_{n,1}](r)
\ls \|\hat \sigma\|_{\U^1_{\alpha+2,\kappa+2}},\\
\displaystyle \M^0_{\alpha+2,\kappa+2}[\Q_{n,2}](r)
\ls \frac{R_*^{\rho_*-\alpha}}{\rho_*-\alpha}\|\hat \sigma\|_{\U^1_{\alpha+2,\kappa+2}}
\end{array}
\end{equation}
for $r\geq1$. By \eqref{P1P2} and \eqref{Q1Q2}, we obtain
\begin{align}
\label{y2}
&\M^0_{\alpha+2,\kappa+2}[y_{\pm2}](r) \\
& \ls
R_*^{\rho_*} \left(\frac1{\alpha(\rho_*-\alpha)} \|\hat w\|^2_{\U^1_{\alpha+2, \kappa+2}}+\nu_*\|\hat w\|_{\U^1_{\alpha+2, \kappa+2}} +\frac{1}{\rho_*-\alpha}\|\hat \sigma\|_{\U^1_{\alpha+2,\kappa+2}}\right) \nonumber
\end{align}
for $r\geq1$. 
On the other hand, in a similar way to Step 2 with regard to the form of $\dr y_n$, we also have
\begin{align}
\label{y22}
&\M^1_{\alpha+2,\kappa+2}[y_{\pm2}](r) \\
& \ls
R_*^{\rho_*} \left(\frac1{\alpha(\rho_*-\alpha)} \|\hat w\|^2_{\U^1_{\alpha+2, \kappa+2}}+\nu_*\|\hat w\|_{\U^1_{\alpha+2, \kappa+2}} +\frac{1}{\rho_*-\alpha}\|\hat \sigma\|_{\U^1_{\alpha+2,\kappa+2}}\right) \ r\geq 1.\nonumber
\end{align}
\vspace{5pt}

{\it Step 3-2: The case $n=\pm2$ and $0\leq r\leq 1$.}\hspace{5pt}
It is sufficient to reconsider $P_{n,1}$ and $Q_{n,1}$. 
Applying a similar discussion to Step 2, we also obtain the estimates for $0\leq r\leq1$ such that
\begin{equation}
\label{y23}
\M^0_{\alpha+2,\kappa+2}[y_{\pm2}](r) \ls
 \left(\frac1{\alpha} \|\hat w\|^2_{\U^1_{\alpha+2, \kappa+2}}+\nu_*\|\hat w\|_{\U^1_{\alpha+2, \kappa+2}} \right)r +\|\hat \sigma\|_{\U^1_{\alpha+2,\kappa+2}}
\end{equation}
and
\begin{equation}
\label{y24}
\M^1_{\alpha+2,\kappa+2}[y_{\pm2}](r) \ls
\left(\frac1{\alpha} \|\hat w\|^2_{\U^1_{\alpha+2, \kappa+2}}+\nu_*\|\hat w\|_{\U^1_{\alpha+2, \kappa+2}} \right)+\|\hat \sigma\|_{\U^1_{\alpha+2,\kappa+2}} .
\end{equation}

\eqref{y0}, \eqref{dry0}, \eqref{y0}-\eqref{ynot22}, and \eqref{y2}-\eqref{y24}
 completes the proof of Proposition \ref{prop:4}. \qed

\renewcommand{\arraystretch}{1}

\subsection{Proof of Theorem \ref{thm:1}}
Finally, let us prove our main theorem.
In what follows, let $\kappa>1$.

{\it Step 1: Existence and regularity of solutions of \eqref{VSp}.}

By Proposition \ref{prop:4}, there is a constant $C_0>0$ such that
\begin{align*}
&\|\S(\hat w, \hat\sigma)\|_{\U^1_{\alpha+2,\kappa+2}}\\
& \leq
C_0 R_*^{\rho_*} \left(\frac1{\alpha(\rho_*-\alpha)} \|\hat w\|^2_{\U^1_{\alpha+2, \kappa+2}}+\nu_*\|\hat w\|_{\U^1_{\alpha+2, \kappa+2}} +\frac{1}{\rho_*-\alpha}\|\hat \sigma\|_{\U^1_{\alpha+2,\kappa+2}}\right) \\
& =: K_1\|\hat w\|^2_{\U^1_{\alpha+2, \kappa+2}}+K_2\nu_*\|\hat w\|_{\U^1_{\alpha+2, \kappa+2}}
+K_3\|\hat \sigma\|_{\U^1_{\alpha+2,\kappa+2}}
\end{align*}
for every $0<\alpha<\min\{1/2, \rho_*\}$, $\hat w\in \U^0_{\alpha+2, \kappa+2}$, and $\hat\sigma\in \U^1_{\alpha+2,\kappa+2}$. 
We should note that $K_1$ and $K_3$ are constants dependent on $R_*, \rho_*$, and $\alpha$, while $K_2$ is a constant dependent only on $R_*$ and $\rho_*$. 
Then in advance, we set $\phi_*$ so that
\begin{equation*}
R_*^{\rho_*}\nu_*< \frac{1}{C_0}=:\delta,
\end{equation*}
in order to see $K_2\nu_*<1$. Moreover, let $\varepsilon>0$ be such as
\begin{equation*}
\varepsilon:=\frac{(1-K_2\nu_*)^2}{4K_1K_3}
\end{equation*}
so that
\[
(1-K_2\nu_*)^2-4K_1K_3\|\hat \sigma\|_{\U^1_{\alpha+2,\kappa+2}}>0
\]
for every $\hat\sigma\in \U^1_{\alpha+2,\kappa+2}$ such that $\|\hat \sigma\|_{\U^1_{\alpha+2,\kappa+2}}<\varepsilon$.

We now take $\varphi \in C^2(\R^2;\R)$ such that 
\begin{equation}
\label{F}
\|\hat \varphi\|_{\U^1_{\alpha+2,\kappa+2}}<\varepsilon.
\end{equation}
We note here that by continuity of $\varphi$ at $r=0$, it holds automatically that $\varphi_n(0)=0$ for $n\in\Z\bs\{0\}$. 
Therefore, $\tilde\varphi$ automatically satisfies the condition of $\U^1_{\alpha+2,\kappa+2}$ at $r=0$.
Then we define the approximative sequence $(\hat w^{(j)})_{j\in\N}$ to the solution of the equation
\begin{equation}
\label{Eq1}
\hat w= \S(\hat w, \hat\varphi)
\end{equation}
as
\begin{equation*}
\left\{
\begin{array}{ll}
\hat w^{(1)}:= \hat\Phi(\hat\varphi),\\
\hat w^{(j)}:= \S(\hat w^{(j-1)}, \hat\varphi),\hspace{10pt}j\geq2,
\end{array}
\right.
\end{equation*}
where $\hat\Phi(\hat\varphi)=(\Phi_n(\hat\varphi))_{n\in\Z}$ is such that
\begin{equation*}
\Phi_n(\hat\varphi):=
\left\{
\begin{array}{lll}
\varphi_n, &n\in\Z\backslash\{\pm2\},\\
\Q_{n,1}[\varphi_n], &n=\pm2,\ 0\leq r\leq R_*,\\
\Q_{n,2}[\varphi_n], &n=\pm2,\ r>R_*.
\end{array}
\right.
\end{equation*}
Let $M>0$ be as
\begin{equation*}
M:=\frac{1-K_2\nu_*-\sqrt{(1-K_2\nu_*)^2-4K_1K_3\|\hat\varphi\|_{\U^1_{\alpha+2,\kappa+2}}}}{2K_1}.
\end{equation*}
Then we see from the estimate of $\hat\Phi(\hat\varphi)$ in the proof of Proposition \ref{prop:4} (see \eqref{Q1Q2}) that
\begin{equation*}
\|\hat w^{(1)}\|_{\U^1_{\alpha+2,\kappa+2}}
\leq K_3\|\hat\varphi\|_{\U^1_{\alpha+2,\kappa+2}} \leq M,
\end{equation*}
and if $\|\hat w^{(j)}\|_{\U^1_{\alpha+2,\kappa+2}}\leq M$ for some $j\in\N\cup\{0\}$, then
\begin{equation*}
\|\hat w^{(j+1)}\|_{\U^1_{\alpha+2,\kappa+2}}
\leq K_1M^2+K_2\nu_*M+K_3\|\hat \sigma\|_{\U^1_{\alpha+2,\kappa+2}}
=M.
\end{equation*}
Therefore, by induction, we see that the sequence $(\|\hat w^{(j)}\|_{\U^1_{\alpha+2,\kappa+2}})_{j\in\N}$ is uniformly bounded by $M$.
Moreover, since $\L$, $I^T_z$ and $J^t_z$ are linear functionals, and since it holds
\begin{equation*}
\cN(f_1,f_2)-\cN(g_1,g_2)= \cN(f_1, f_2-g_2)+\cN(f_1-g_1, g_2)
\end{equation*}
for every bilinear form $\cN$ appearing in $\G^0_n$, we inductively see
\begin{align*}
\|\hat w^{(j+1)}-\hat w^{(j)}\|_{\U^1_{\alpha+2,\kappa+2}}
&\leq \|\S(\hat w^{(j)}, \hat\varphi)-\S(\hat w^{(j-1)}, \hat\varphi)\|_{\U^1_{\alpha+2,\kappa+2}}\\
&\leq (2K_1M+K_2\nu_*)\|\hat w^{(j)}-\hat w^{(j-1)}\|_{\U^1_{\alpha+2,\kappa+2}}\\
&\leq (2K_1M+K_2\nu_*)^{j-1}\|\S(\hat\Phi(\hat\varphi), \hat\varphi)-\hat\Phi(\hat\varphi)\|_{\U^1_{\alpha+2,\kappa+2}}\\
&\leq (2K_1M+K_2\nu_*)^j M
\end{align*}
by reviewing the proof of Proposition \ref{prop:4}. 
Since $2K_1M+K_2\nu_*<1$, we see
\begin{equation*}
\sum_{j=1}^\infty\|\hat w^{(j+1)}-\hat w^{(j)}\|_{\U^1_{\alpha+2,\kappa+2}}<\infty,
\end{equation*}
which and the completeness of $\U^1_{\alpha+2,\kappa+2}$ yield that $\hat w^{(j)}$ converges to some $\hat w^{\infty} \in \U^1_{\alpha+2,\kappa+2}$ under the condition \eqref{F}. Since
\begin{equation*}
\|\S(\hat w^{\infty} , \hat\varphi)-\S(\hat w^{(j)} , \hat\varphi)\|_{\U^1_{\alpha+2,\kappa+2}}
\leq (2K_1M-K_2\nu_*)\|\hat w^{\infty}-\hat w^{(j)}\|_{\U^1_{\alpha+2,\kappa+2}}\to 0
\end{equation*}
as $j\to\infty$, we see that $\hat w^{\infty}$ is a solution of \eqref{Eq1}.

In what follows, we fix a pair $(\hat w, \hat \varphi) \in \U^1_{\alpha+2,\kappa+2}\times \U^1_{\alpha+2,\kappa+2}$ of solutions obtained as above. 
Then this pair satisfies
\renewcommand{\arraystretch}{1.5}
\begin{equation*}
w_n(r):=
\left\{
\begin{array}{llll}
 \displaystyle    -I^\infty_{|n|}[\G^*_n](r)- J^0_{|n|}[\G^*_n](r)+\varphi_n(r),  &r\geq0,\ n\in\Z\backslash\{0,\pm2\},\\
\displaystyle  w_{n,1}(r):=\P_{n,1}(r)+ \Q_{n,1}[\varphi_n](r)
&0\leq r\leq R_*,\ n=\pm2,\\
\displaystyle w_{n,2}(r):=\P_{n,2}(r)+ \Q_{n,2}[\varphi_n](r)
&r> R_*,\ n=\pm2,\\
\displaystyle   -\int^\infty_r\frac1s\int_0^s t\G^*_0(t)dtds+\varphi_n(r),  &r\geq0, \ n=0,
  \end{array}
\right.
\end{equation*}
together with $\hat\gamma:=\L(\hat w) \in \V^2_{\alpha, \kappa+4}$.
Since $\hat w\in \U^1_{\alpha+2,\kappa+2}$, $\hat\varphi\in \hat C^2$, and since $\G^*_n=\H_n$ is finite at $r=0$ (see \eqref{H} in Proposition \ref{prop:2}), 
we see that each of $\Delta_{r,n}w_n$ $(|n|\neq2)$, $\Delta_{r, n}w_{n,1}$ and $\Delta_{r, \zeta_{n}}w_{n,2}$ $(|n|=2)$ is well-defined and continuous in each domain, and is expressed as
\begin{alignat}{2}
&\Delta_{r,n}w_n(r)= \G^*_n(r)+ \Delta_{r,n}\varphi_n(r)&,\hspace{10pt}&0\leq r<\infty,\ n\in\Z\bs\{\pm2\},
\label{wnot2}\\
&\Delta_{r, n}w_{n,1}(r)=\G^*_n(r)+ \Delta_{r,n}\varphi_n(r)&,\hspace{10pt}&0\leq r\leq R_*,\ |n|=2,
\label{wn1}\\
&\Delta_{r, \zeta_{n}}w_{n,2}(r)=\G^0_n(r)+ \Delta_{r,n}\varphi_n(r)&,\hspace{10pt}&R_*< r<\infty,\ |n|=2.
\label{wn2}
\end{alignat}
From \eqref{wn1}, \eqref{wn2}, and the connecting properties \eqref{match1} and \eqref{match2}, we see
\begin{equation*}
\dr^l w_{n,1}(R_*)=\lim_{r\to R_*-0}\dr^l w_{n,2}(r),\hspace{5pt}l=0,1,2,
\end{equation*}
and hence each of $\Delta_{r,n}w_{\pm2}$ is also well-defined and continuous in $[0,\infty)$, and satisfies \eqref{wnot2}.
Therefore, we see that $(\hat\gamma, \hat w)=(\L(\hat w), \hat w)$ are strong solutions of the system \eqref{VSp}.\\

{\it Step 3: Existence of a solution of \eqref{NS} and its decay property.}\hspace{5pt}
Using the above solutions $(\hat\gamma, \hat w)$ of \eqref{VSp}, we set
\begin{equation*}
\gamma(r,\theta):=\sum_{n\in\Z}\gamma_n(r)e^{in\theta}, \hspace{10pt}
w(r,\theta):=\sum_{n\in\Z}w_n(r)e^{in\theta}
\end{equation*}
for every $(r,\theta)\in[0,\infty)\times [0,2\pi)$. By the definition of $\V^2_{\alpha, \kappa+4}$ and 
$\U^1_{\alpha+2,\kappa+2}$, and the discussion in Step 2, 
we see $\gamma \in C^2(\R^2; \R)$ and $w\in C^2(\R^2;\R)$.
Moreover, by summing up Fourier modes of the system \eqref{VSp}, we also see that $(\psi, \omega)$ defined by
\begin{equation*}
\psi(r,\theta):=\psi_*(r)+\gamma(r,\theta),\hspace{5pt}
\omega(r,\theta):=\omega_*(r)+w(r,\theta).
\end{equation*}
are strong solutions of the vorticity-streamfunction system \eqref{VS2}, i.e., \eqref{VS} for $\phi=\phi_*+\varphi$.
We note here that by the elliptic regularity of the Poisson equation, we see that $\psi \in C^3(\R^2, \R)$.

\renewcommand{\arraystretch}{1}
Now let $u:=\nabla^\perp\psi$. 
Then we see $u\in C^2(\R^2, \R)$. 
On the other hand, $\dr\gamma_0$ is actually expressed as
\begin{align*}
\dr\gamma_0(r)
&=-\frac1r\int_0^r\left[-\D_0(s)+2s\int_s^\infty \frac{\D_0(t)}{t^2}dt +s\varphi_0(s)\right]ds\\
&=-r\int_r^\infty\frac{\D_0(s)}{s^2}ds -\frac1r\int_0^r s\varphi_0(s)ds.
\end{align*}
We see from \eqref{dn0at0} that the modulus of this first term is bounded by $(1+r)^{-(\alpha+1)}$. 
Hence together with $\tilde\gamma\in\tilde{V}^2_{\alpha, \kappa+4}$, 
we obtain the decay property \eqref{decay} of $u$ by using the formula $\nabla^\perp(\phi_*+\varphi_0)=-\dr(\phi_*+\varphi_0)e_\theta$. 
Furthermore, this $u$ and $\omega$ satisfy
\begin{equation*}
\left\{
  \begin{array}{ll}
    \Delta \omega=u\cdot\nabla \omega+ \Delta (\phi_*+\varphi), \hspace{5pt}  &x\in \R^2, \\
    \nabla\cdot u=0 \hspace{5pt}  &x\in \R^2,
  \end{array}
\right.
\end{equation*}
and the decay properties
\begin{equation*}
|u(r,\theta)|\ls(1+r)^{-1},\hspace{5pt}|\nabla u(r,\theta)|\ls(1+r)^{-2}, \hspace{5pt}
|\nabla\omega(r,\theta)|\ls (1+r)^{-3}.
\end{equation*}
Hence using a similar method to Hillairet-Wittwer \cite[Section 3]{HW2013}, 
we see that $u$ becomes a classical solution of the Navier-Stokes system \eqref{NS} for $F=\nabla^\perp (\phi_*+\varphi)$ together with some pressure $p$.
\qed

\section*{Acknowledgement}
The work of Y.~M. is partially supported by JSPS KAKENHI (Grant Number:  19H05597, 20H00118, 20K03698).
The work of H.~T. is partially  supported by Grant-in-Aid for JSPS Research Fellow (Grant Number: 21J00379).
The authors would like to express their sincere thanks to the referees for their valuable comments and advices.

\renewcommand{\refname}{References}

\end{document}